\documentclass{amsart}

\usepackage{amsmath,amsthm,amssymb,amscd,enumerate,mathtools,enumitem,tikz-cd,bm}
\usepackage{amsfonts}
\usepackage{rotating}
\usepackage{euscript}
\usepackage{pst-node}
\usepackage{epsfig,verbatim}
\newcommand{\href}[2]{#2}

\def\be{\begin{equation}}
\def\ee{\end{equation}}

\def\C{{\mathbb C}} 

\def\N{{\mathbb N}} 
\def\P{{\mathbb P}}

\def\phi{{\varphi}}

\def\deg{{\rm deg\,}} 

\def\dim{{\rm dim\,}}

\def\bp{\begin{proposition}}
\def\ep{\end{proposition}}

\def\bt{\begin{theorem}}
\def\et{\end{theorem}}
\def\br{\begin{remark}}
\def\er{\end{remark}}
\def\be{\begin{equation}}
\def\bee{\begin{equation*}}
\def\l{\label}

\def\ee{\end{equation}}
\def\eee{\end{equation*}}
\def\bl{\begin{lemma}}
\def\el{\end{lemma}}
\def\bc{\begin{corollary}}
\def\ec{\end{corollary}}
\def\pr{\noindent{\it Proof. }}

\def\bd{\begin{definition}}
\def\ed{\end{definition}}

\newtheorem{theorem}{Theorem}[section]
\newtheorem{lemma}[theorem]{Lemma}
\newtheorem{definition}[theorem]{Definition}
\newtheorem{corollary}[theorem]{Corollary}
\newtheorem{proposition}[theorem]{Proposition}
\newtheorem{problem}[theorem]{Problem}

\theoremstyle{definition}
\newtheorem{example}[theorem]{Example}
\theoremstyle{definition}
\newtheorem{remark}[theorem]{Remark}
\def\bpr{\begin{problem}}
\def\epr{\end{problem}}

\newcommand{\Cent}{\mathcal C}

\DeclareMathOperator{\per}{per}

\mathtoolsset{showonlyrefs}

\title[Bounds on periods of curves]{Bounds on periods of curves for split\\ rational maps}
\author{}
\date{}

\begin{document}

\author[F. Pakovich]{Fedor Pakovich}
\thanks{
This research was supported by ISF Grant No. 1092/22} 
\address{Department of Mathematics, 
Ben Gurion University of the Negev, P.O.B. 653, Beer Sheva,  8410501, Israel}
\email{pakovich@math.bgu.ac.il}

\begin{abstract}
Let $A_1$ and $A_2$ be fixed complex rational maps of degree
at least two, and let $(A_1,A_2)$ denote the endomorphism of $\P^1\times\P^1$ of the form 
\linebreak $(x,y)\mapsto(A_1(x),A_2(y))$. We show that the periods of irreducible
$(A_1,A_2)$-periodic curves with dominant coordinate projections are
bounded, unless $A_1$ and $A_2$ belong to the exceptional power,
Chebyshev, or Latt\`es families. 
We also extend this result to periodic subvarieties  for endomorphisms of 
$(\P^1)^m$.

\end{abstract}

\maketitle

\section{Introduction}
Let $A_1$ and $A_2$ be complex rational maps of degree at least two, and let $(A_1,A_2)$ denote the  
endomorphism of $\P^1\times\P^1$ of the form 
\be \l{end} \qquad (x,y)\longmapsto(A_1(x),A_2(y)).
\ee 
An irreducible algebraic curve $C$ in $\P^1\times\P^1$ is called \emph{$(A_1,A_2)$-periodic} if
$$(A_1,A_2)^{\circ m}(C)=C$$ for some $m\geq1$, and the least such integer is called the \emph{period} of $C$,
denoted by $\per_{(A_1,A_2)}(C)$. Throughout the paper we are interested
in curves whose two coordinate projections are dominant, or equivalently,
curves that are neither vertical nor horizontal.

When $A_1$ and $A_2$ are polynomials, Medvedev and Scanlon
\cite{MS} classified invariant curves for endomorphisms
of the form \eqref{end} using Ritt's theory of polynomial
decompositions \cite{r1}. Their results also describe
periodic curves by passage to iterates. For rational
maps, a description of invariant, periodic, and preperiodic curves
was obtained in \cite{PakInv}, based on the theory of semiconjugate
rational functions developed in \cite{semi,rec,fin,PakLat}.

A natural question is whether the periods of $(A_1,A_2)$-periodic
curves can be bounded and, if so, in terms of which data. Without further assumptions on $A_1$ and $A_2$, such bounds need
not exist, owing to the presence of special maps. 
Recall that a rational map of degree $d\geq2$ is called \emph{special} if it is
a Latt\`es map or is M\"obius conjugate to $z^{\pm d}$ or $\pm T_d$,
where $T_d$ is the Chebyshev polynomial.  

For a fixed special map
$A$, already the diagonal map $(A,A)$ admits  irreducible periodic curves
with dominant coordinate projections and arbitrarily large periods. For power maps and Chebyshev polynomials, such examples can be found
in \cite[p.~214]{GN} and \cite[Example~5.3]{PakPoly}, respectively. For Latt\`es maps, we construct them below. 

Note that the boundedness question is meaningful only for curves with dominant
coordinate projections, since curves of the form $\{x_0\}\times\P^1$, where $x_0$ is a periodic
point of $A_1$, have the same period as $x_0$, which can be arbitrarily large.

For non-special  polynomials $A_1$ and $A_2$, the periods of periodic
curves with dominant coordinate projections admit bounds depending
only on the degrees of $A_1$ and $A_2$. The existence of such a curve forces
$\deg A_1=\deg A_2$, and writing $d$ for this common degree, Ghioca and
Nguyen established the explicit bound $d^{2d^4}/2$
in \cite[Theorem~1.1]{GN}.
In the recent work \cite[Theorem~1.2]{PakPoly}, this bound was improved to
\be \l{bou}  
\operatorname{lcm}\bigl(|\Sigma(A_1)|,|\Sigma(A_2)|\bigr),
\ee
where $\Sigma(A)$ is the symmetry group of the Julia set $J(A)$,
consisting of all affine transformations $\mu(z)=az+b$, $a\ne0$,
such that $\mu(J(A))=J(A)$. 
Since $\vert \Sigma(A_i)\vert \leq d,$ this result gives the bound $d^
2$
in terms of the common
degree. The bound \eqref{bou} was also shown to be optimal:
it is attained for polynomials of arbitrarily large degree with
arbitrarily large symmetry groups.

In the rational setting, analogous bounds depending only on the degrees
of the maps are not currently known. However, some related results
are available. For example, Bell, Ghioca, and Tucker \cite{BGT} bounded
the periods of periodic curves under split rational
maps over number fields in terms of the degrees,
the field of definition, and a chosen place of good
reduction of the maps. Further, \cite[Theorem~1.2]{DM} implies that, for non-special
rational maps $A_1,A_2$, the periods of irreducible periodic
curves of a given bidegree with dominant coordinate projections
are bounded in terms of this bidegree and the degrees of the maps. On the other hand, it was shown
in \cite[Theorem~1.2]{PakGen} that, for every $d\geq2$, there exists
a nonempty Zariski open subset $U\subset\operatorname{Rat}_d$ such
that, for any $A_1,A_2\in U$, every irreducible
$(A_1,A_2)$-periodic curve with dominant coordinate projections
is actually invariant.

The above results leave open the following  question:
can a fixed pair of non-special rational maps $A_1,A_2$ admit
irreducible $(A_1,A_2)$-periodic curves with dominant coordinate
projections and arbitrarily large periods? 
In this note, we answer this question in the negative by combining
the classification of invariant curves in \cite{PakInv} with a recent
result of Beaumont \cite[Corollary~1.4]{Beaumont}, which states that
all rational maps commuting with some iterate of a fixed non-special
rational map already commute with one common iterate.

Specifically, we prove the following result.

\begin{theorem}\label{thm:main} 
For every pair $A_1,A_2$ of non-special complex rational maps of degree at least two, 
there exists an integer $N\geq 1$ such that every irreducible
$(A_1,A_2)$-periodic curve $C\subset\P^1\times\P^1$ with dominant coordinate
projections satisfies
\begin{equation}\label{eq:main}
 (A_1^{\circ N},A_2^{\circ N})(C)=C.
\end{equation}
In particular, $\per_{(A_1,A_2)}(C)$ divides $N$. The integer $N$ depends on $A_1$ and $A_2$ but is independent
of $C$ and of its bidegree.
\end{theorem}

The results of \cite{BGT,DM} cited above also concern
subvarieties of higher dimension. We obtain an analogous extension
of Theorem~\ref{thm:main} using the results of \cite{MS},
which reduce the description of invariant subvarieties
of $(\P^1)^m$, $m\geq2$, for split maps  to the description
of invariant curves in $\P^1\times\P^1$.
Note, however, that neither the proof of Theorem~\ref{thm:main} nor that
of its generalization yields an effective bound on the periods.

This note is organized as follows. Section~2 recalls the results on semiconjugate rational maps, periodic curves, and centralizers needed in the proof.
Section~3 proves Theorem~\ref{thm:main}. Finally, Section~4 extends the result
to periodic subvarieties of $(\P^1)^m$.

\section{Preliminaries}\label{sec:prelim}

\subsection{Semiconjugate rational maps}
Below, all maps are assumed to be nonconstant. 
We start by recalling the notion of a generalized Latt\`es map,
introduced in \cite{PakLat}, which is very useful when describing rational semiconjugacies.

An orbifold on $\P^1$
is specified by a {\it ramification function}
$\nu\colon\P^1\to\N$ equal to $1$ outside a finite set.
We consider only good orbifolds: the set
$\{z\in\P^1:\nu(z)>1\}$ must not consist of one point
or of two points with unequal values of $\nu$. A rational
map $A$ of degree at least two is called a \emph{generalized Latt\`es map}
if there is such a function $\nu\not\equiv1$ satisfying
\[
 \nu(A(z))=\nu(z)\gcd\bigl(\deg_z A,\nu(A(z))\bigr),
 \qquad z\in\P^1,
\]
where $\deg_z A$ is the local degree of $A$ at $z$. Every special
map is a generalized Latt\`es map, but the latter class is much larger. 
A rational map $A$ is a generalized Latt\`es map if and only if
some positive iterate of $A$ is one (see \cite[Theorem~2.7]{PakInv}
and the discussion following it).

We will use the following consequence of
\cite[Proposition~3.3]{PakInv}.

\begin{theorem}\label{thm:reverse-semi}
Let $A,B$ be rational maps of degree at least two, and let $X$ be
a rational map that makes the diagram
\[
\begin{CD}
\P^1 @>{B}>> \P^1\\
@V{X}VV @VV{X}V\\
\P^1 @>{A}>> \P^1
\end{CD}
\]
commute.
Suppose that $A$ is not a generalized Latt\`es map.
Then there exists a rational map $Y$ such that the diagram
\[
\begin{CD}
\P^1 @>{A}>> \P^1\\
@V{Y}VV @VV{Y}V\\
\P^1 @>{B}>> \P^1
\end{CD}
\]
commutes.
\end{theorem}

The following result allows us to lift a non-special rational
map through a semiconjugacy to a map which is not a generalized
Latt\`es map (see \cite[Theorem~4.14]{PakInv} and the first
paragraph on p.~271 of the same paper).

\begin{theorem}\label{thm:cover}
Let $B$ be a non-special rational map of degree
at least two. Then there exist rational functions $F$ and
$\theta$ such that $F$ is not a generalized Latt\`es map, $\theta$ is a Galois covering,
 and the diagram
\begin{equation}\label{eq:cover}
\begin{CD}
\P^1 @>{F}>> \P^1\\
@V{\theta}VV @VV{\theta}V\\
\P^1 @>{B}>> \P^1
\end{CD}
\end{equation}
commutes.
\end{theorem}

\subsection{Invariant and periodic curves}\label{subsec:curves}

Recall that, for non-special rational maps $A_1,A_2$ of degree
at least two, an irreducible curve $C\subset\P^1\times\P^1$
with dominant coordinate projections is $(A_1,A_2)$-invariant
if and only if there exist rational maps $X_1,X_2,$ and $B$ such that
the diagram
\[
\begin{CD}
(\P^1)^2 @>{(B,B)}>> (\P^1)^2\\
@V{(X_1,X_2)}VV @VV{(X_1,X_2)}V\\
(\P^1)^2 @>{(A_1,A_2)}>> (\P^1)^2
\end{CD}
\]
commutes and
\[
t\longmapsto(X_1(t),X_2(t))
\]
is a parametrization of $C$ (see \cite[Theorem~4.1]{PakInv}). Note that the existence of such a diagram implies
$\deg A_1=\deg A_2$.

The following statement is contained in
\cite[Theorem~4.10]{PakInv}. We include its proof because
we will use the same argument for curves in $(\P^1)^m$.

\begin{proposition}\label{cor:diagonal-param}
Let $A$ be a rational map of degree at least two that
is not a generalized Latt\`es map. Then every irreducible
$(A,A)$-periodic curve $C$ with dominant coordinate
projections has a parametrization
\[
 C=\{(U_1(t),U_2(t)):t\in\P^1\},
\]
where each of the maps $U_1,U_2$ commutes with
a positive iterate of $A$.
\end{proposition}

\begin{proof}
Let $s$ be the period of $C$. Then there exist rational maps $V_1,V_2,$ and $H$
such that $t\mapsto(V_1(t),V_2(t))$ parametrizes $C$
and the lower square of the  diagram 
\[
\begin{CD}
(\P^1)^2 @>{(A^{\circ s},A^{\circ s})}>> (\P^1)^2\\
@V{(Y,Y)}VV @VV{(Y,Y)}V\\
(\P^1)^2 @>{(H,H)}>> (\P^1)^2\\
@V{(V_1,V_2)}VV @VV{(V_1,V_2)}V\\
(\P^1)^2 @>{(A^{\circ s},A^{\circ s})}>> (\P^1)^2 
\end{CD}
\]
commutes.
Since $A^{\circ s}$ is not a generalized Latt\`es map,
Theorem~\ref{thm:reverse-semi}, applied to the first
coordinate semiconjugacy, gives a rational map $Y$
making the upper square commute.  
Thus the functions $U_i=V_i\circ Y$, $i=1,2$, commute
with $A^{\circ s}$. Since $Y$ is surjective, they give
the required parametrization of $C$. 
\end{proof}

The following simple result  
allows us, to a certain extent, to reduce the study of
$(A_1,A_2)$-periodic curves to the case $A_1=A_2=A$ (see \cite[Lemma~4.4]{PakInv}). 

\begin{lemma}\label{lem:lift}
Let $A_1,A_2,B_1,B_2,X_1,X_2$ be nonconstant rational
functions such that the diagram
\[
\begin{CD}
(\P^1)^2 @>{(B_1,B_2)}>> (\P^1)^2\\
@V{(X_1,X_2)}VV @VV{(X_1,X_2)}V\\
(\P^1)^2 @>{(A_1,A_2)}>> (\P^1)^2
\end{CD}
\]
commutes. Then for every irreducible $(A_1,A_2)$-periodic
 curve $C$ there exists an irreducible
$(B_1,B_2)$-periodic curve $D$ such that
$
 C=(X_1,X_2)(D).
$
\end{lemma}

\subsection{Centralizers of iterates}

For a rational map $A$ of degree at least two, set
\[
 \Cent(A)=\{U\in\C(z)\setminus\C:U\circ A=A\circ U\}
\] 
and 
\[
 \qquad \Cent_\infty(A)=\bigcup_{k\geq1}\Cent(A^{\circ k}).
\]
It is easy to see that both sets are semigroups under composition. Furthermore,
the classical theorem of Ritt \cite{RittPerm} implies that,
if $A$ is non-special, a rational map $U$ of degree at least
two belongs to $\Cent_\infty(A)$ if and only if
$U^{\circ k}=A^{\circ \ell}$ for some $k,\ell\geq1$.  

It was proved in \cite[Theorem~1.2]{fin} that, for non-special
$A$, there exist finitely many rational maps $U_1,\ldots,U_s$
such that every element of $\Cent(A)$ has the form
$U_j\circ A^{\circ k}$ with $k\geq0$. Further results on the
structure of $\Cent(A)$ were obtained in \cite{PakComm}.

In \cite[Conjecture~3.1 and Conjecture~3.2]{PakConj}, we conjectured the
analogous finiteness statement for $\Cent_\infty(A)$ and 
Beaumont proved this conjecture,  
establishing the following result (see  \cite[Theorem~1.3]{Beaumont} and
\cite[Corollary~1.4]{Beaumont}).

\begin{theorem}\label{thm:beaumont}
Let $A$ be a rational map of degree at least two that
is not M\"obius conjugate to $z^{\pm\deg A}$.
Then there exists an integer $n(A)\geq1$ such that
\begin{equation}\label{eq:centralizer}
 \Cent_\infty(A)=\Cent(A^{\circ n(A)}).
\end{equation}
\end{theorem}

\section{Proof of Theorem~\ref{thm:main}}
Assume that there exists at least one irreducible $(A_1,A_2)$-periodic
curve with dominant coordinate projections. Fix such a curve $C_0$
and denote its period by $r$. Choose rational functions $X_1,X_2,B$
such that the map $t\mapsto (X_1(t),X_2(t))$ parametrizes $C_0$
and the  diagram 
\[
\begin{CD}
(\P^1)^2 @>{(B,B)}>> (\P^1)^2\\
@V{(X_1,X_2)}VV @VV{(X_1,X_2)}V\\
(\P^1)^2 @>{(A_1^{\circ r},A_2^{\circ r})}>> (\P^1)^2
\end{CD}
\]
commutes. 
Since the map $A_1$ is non-special, its iterate $A_1^{\circ r}$
is also non-special (see, e.g., \cite[Lemma~2.12]{fin}). Thus, since $B$ is semiconjugate to $A_1^{\circ r}$, the map $B$
is also non-special by \cite[Corollary~4.3]{Milnor}
(this also follows from \cite[Theorem~1.14]{JX}, using
the relation between multipliers of periodic points
of semiconjugate maps).

Applying Theorem~\ref{thm:cover} to $B$, we obtain
rational functions $\theta$ and $F$ such that $F$ is not
a generalized Latt\`es map and the diagram
\[
\begin{CD}
(\P^1)^2 @>{(F,F)}>> (\P^1)^2\\
@V{(\theta,\theta)}VV @VV{(\theta,\theta)}V\\
(\P^1)^2 @>{(B,B)}>> (\P^1)^2\\
@V{(X_1,X_2)}VV @VV{(X_1,X_2)}V\\
(\P^1)^2 @>{(A_1^{\circ r},A_2^{\circ r})}>> (\P^1)^2
\end{CD}
\]
commutes. Note that all the maps in this diagram are fixed for the pair $(A_1,A_2)$ once a curve $C_0$ is chosen. 

By Theorem~\ref{thm:beaumont}, there exists an integer $n\geq1$
such that every rational function commuting with an iterate
of $F$ commutes with $F^{\circ n}$. We claim that $N=rn$
satisfies the conclusion of the theorem.

Indeed, let $C$ be any irreducible $(A_1,A_2)$-periodic curve
with dominant coordinate projections. Since $C$ is also
$(A_1^{\circ r},A_2^{\circ r})$-periodic,
Lemma~\ref{lem:lift}, applied to the composite diagram,
gives an irreducible $(F,F)$-periodic curve $D$ such that
\[
 C=(X_1\circ\theta,X_2\circ\theta)(D).
\]
Moreover, since the coordinate projections of $C$ are dominant, 
the same holds for $D$.

By Proposition~\ref{cor:diagonal-param}, the curve $D$ has
a parametrization $t\mapsto(U_1(t),U_2(t))$, where each
of the maps $U_1,U_2$ commutes with an iterate of $F$.
Consequently, both commute with $F^{\circ n}$, and hence
\[
\begin{aligned}
 (F^{\circ n},F^{\circ n})(D)
 &=\{(U_1(F^{\circ n}(t)),
       U_2(F^{\circ n}(t))):t\in\P^1\}=D.
\end{aligned}
\]
Setting $$\pi=(X_1\circ\theta,X_2\circ\theta),$$ we therefore obtain
$$
\begin{aligned}
 (A_1^{\circ rn},A_2^{\circ rn})(C)
 &=\pi\bigl((F^{\circ n},F^{\circ n})(D)\bigr)=\pi(D)=C.
\end{aligned} \eqno{\Box} 
$$

\begin{example}
Let $A$ be a Latt\`es map. By \cite[Theorem~3.1]{Milnor},
there exist a complex torus $E=\C/\Lambda$, a finite
cyclic group $G$ of rotations of $E$ about the origin,
and an affine map $\ell:E\to E$ commuting with $G$,
such that for the quotient map $$\rho:E\to E/G\cong\P^1$$
the equality 
\[
 A\circ\rho=\rho\circ\ell 
\] holds. 
Write $\ell(P)=aP+b$ and let $\alpha:E\to E$ be its
linear part, $\alpha(P)=aP$.

For $T\in E$, define
\[
 C_T=\{(\rho(P),\rho(P+T)):P\in E\}.
\]
Each $C_T$ is an irreducible curve with dominant coordinate
projections. For every $k\geq1$, we have
\[
 \ell^{\circ k}(P)=a^kP+(1+a+\cdots+a^{k-1})b,
\]
and hence
\[
 \ell^{\circ k}(P+T)
 =\ell^{\circ k}(P)+a^kT
 =\ell^{\circ k}(P)+\alpha^{\circ k}(T).
\]
Using this identity and the surjectivity of $\ell^{\circ k}$,
we obtain
\[
\begin{aligned}
 (A,A)^{\circ k}(C_T)
 &=\{(\rho(\ell^{\circ k}(P)),\rho(\ell^{\circ k}(P+T))):P\in E\}\\
 &=\{(\rho(\ell^{\circ k}(P)),
       \rho(\ell^{\circ k}(P)+\alpha^{\circ k}(T))):P\in E\}\\
 &=\{(\rho(Q),\rho(Q+\alpha^{\circ k}(T))):Q\in E\}\\
 &=C_{\alpha^{\circ k}(T)},\qquad k\geq1.
\end{aligned}
\]
Moreover,
\[
 C_S=C_T \quad\Longleftrightarrow\quad \rho(S)=\rho(T).
\]
Indeed, since every element of $G$ fixes the origin,
the only point of $C_T$ with first coordinate $\rho(0)$
is $(\rho(0),\rho(T))$. Thus $C_S=C_T$ implies
$\rho(S)=\rho(T)$. Conversely, if $\rho(S)=\rho(T)$,
then $S=\varepsilon T$ for some $\varepsilon\in G$,
and the substitution $P\mapsto\varepsilon P$ gives
$C_S=C_T$.

The preceding equalities show that
\[
 (A,A)^{\circ k}(C_T)=C_T
 \quad\Longleftrightarrow\quad
 \rho(\alpha^{\circ k}(T))=\rho(T).
\]
Observe now that $\alpha$,
being multiplication by $a$, commutes with every rotation
in $G$. Hence it induces a rational map $A_0$ on
$E/G\cong\P^1$ satisfying
\[
 A_0\circ\rho=\rho\circ\alpha.
\]
Thus $C_T$ is invariant under $(A,A)^{\circ k}$ if and only if
\[
 A_0^{\circ k}(\rho(T))=\rho(T).
\]
Consequently, $C_T$
is periodic under $(A,A)$ if and only if $\rho(T)$ is
periodic under $A_0$, and their exact periods coincide.

Since $\deg A_0=\deg\alpha\geq2$, the map $A_0$ has
periodic points of arbitrarily large periods. Therefore,
the fixed diagonal map $(A,A)$ admits irreducible periodic
curves with dominant coordinate projections and arbitrarily
large periods.
\end{example}

\section{Periodic subvarieties of $(\P^1)^m$}
\label{sec:higher}

Every dominant endomorphism of $(\P^1)^m$ has the form
\[
 (x_1,\ldots,x_m)\longmapsto
 \bigl(R_1(x_{\sigma(1)}),\ldots,R_m(x_{\sigma(m)})\bigr),
\]
where $\sigma$ is a permutation and the $R_i$ are nonconstant
rational maps. Hence some iterate of this endomorphism
is a split rational map of the form
\[
 \Phi:(x_1,\ldots,x_m)\longmapsto
 \bigl(A_1(x_1),\ldots,A_m(x_m)\bigr),
\]
where $A_1,\ldots,A_m$ are rational maps.
We therefore restrict attention to such maps.

For a nonempty subset $I\subseteq\{1,\ldots,m\}$, write
$(\P^1)^I$ for the product of the coordinate factors indexed
by $I$, and put $\Phi_I=(A_i)_{i\in I}$.
We use the following special case of a result of Bell,
Ghioca, and Tucker \cite[Lemma~4.2]{BGT}, which they deduce
from the classification of invariant subvarieties obtained
by Medvedev and Scanlon \cite{MS}.

\begin{lemma}\label{lem:product-curves}
Let $A_1,\ldots,A_m$ be non-special rational maps of degree
at least two, and put $\Phi=(A_1,\ldots,A_m)$.
For every irreducible $\Phi$-periodic subvariety
$V\subset(\P^1)^m$ with dominant coordinate projections,
there exists a partition
\[
 \{1,\ldots,m\}=I_1\sqcup\cdots\sqcup I_r,
 \qquad r=\dim V,
\]
such that, after grouping the coordinates according to this
partition,
\[
 V=C_1\times\cdots\times C_r,
\]
where each $C_j\subset(\P^1)^{I_j}$ is an irreducible
$\Phi_{I_j}$-periodic curve with dominant coordinate
projections.
\end{lemma}

A modification of the proof of Theorem~\ref{thm:main},
combined with Lemma~\ref{lem:product-curves},
yields the following result.

\begin{theorem}\label{thm:higher-split}
Let $A_1,\ldots,A_m$ be non-special rational maps of degree
at least two,  and put $\Phi=(A_1,\ldots,A_m)$.
Then there exists an integer $N\geq1$, depending only on
$A_1,\ldots,A_m$, such that every irreducible
$\Phi$-periodic subvariety $V\subset(\P^1)^m$ with dominant
coordinate projections satisfies
$
 \Phi^{\circ N}(V)=V.
$
\end{theorem}

\pr
Lemma~\ref{lem:product-curves} reduces the proof of
Theorem~\ref{thm:higher-split} to the case of periodic curves for each
subtuple $(A_i)_{i\in I}$. Thus it suffices to establish the assertion for curves
in $(\P^1)^m$, with $m$ arbitrary. Fix an irreducible $(A_1,\ldots,A_m)$-periodic curve
$C_0\subset(\P^1)^m$ with dominant coordinate projections
and denote its period by $r$.

The arguments of Section~\ref{sec:prelim} used in the proof
of Theorem~\ref{thm:main} extend to $m$ coordinates as follows.
First, the normalization argument in
\cite[Theorem~4.1]{PakInv} applies without change: there exist
rational functions $X_1,\ldots,X_m,B$ such that the diagram
\[
\begin{CD}
(\P^1)^m @>{(B,\ldots,B)}>> (\P^1)^m\\
@V{(X_1,\ldots,X_m)}VV @VV{(X_1,\ldots,X_m)}V\\
(\P^1)^m @>{(A_1^{\circ r},\ldots,A_m^{\circ r})}>>
(\P^1)^m
\end{CD}
\]
commutes and $$t\mapsto(X_1(t),\ldots,X_m(t))$$
parametrizes $C_0$.
Indeed, the restriction of $\Phi^{\circ r}$ to $C_0$
lifts to a self-map of its normalization of degree at least
two, so the Riemann--Hurwitz formula implies $g(C_0)\leq 1$.
If $g(C_0)=1$, then all the maps
$A_1^{\circ r},\ldots,A_m^{\circ r}$ are Latt\`es maps,
contrary to their non-specialness. Thus the normalization
of $C_0$ is $\P^1$.

Lemma~\ref{lem:lift} also extends to $m$ coordinates:
an iterate preserving the periodic curve acts on the finite
set of irreducible components of its inverse image,
and a periodic component maps onto that curve. Finally, the
proof of Proposition~\ref{cor:diagonal-param} applies to
$(A,\ldots,A)$-periodic curves in $(\P^1)^m$: the map $Y$
is obtained from the first coordinate semiconjugacy and
then composed with each of $V_1,\ldots,V_m$.

As in the proof of Theorem~\ref{thm:main}, the map $B$
is non-special. Applying Theorem~\ref{thm:cover} to $B$,
we obtain rational maps $\theta$ and $F$ such that $F$
is not a generalized Latt\`es map and the diagram
\[
\begin{CD}
(\P^1)^m @>{(F,\ldots,F)}>> (\P^1)^m\\
@V{(\theta,\ldots,\theta)}VV @VV{(\theta,\ldots,\theta)}V\\
(\P^1)^m @>{(B,\ldots,B)}>> (\P^1)^m\\
@V{(X_1,\ldots,X_m)}VV @VV{(X_1,\ldots,X_m)}V\\
(\P^1)^m @>{(A_1^{\circ r},\ldots,A_m^{\circ r})}>>
(\P^1)^m
\end{CD}
\]
commutes. Put
\[
 \pi=(X_1\circ\theta,\ldots,X_m\circ\theta),
\]
and choose $n=n(F)$ as in Theorem~\ref{thm:beaumont}.

Let $C$ be any irreducible $\Phi$-periodic curve with
dominant coordinate projections. Since $C$ is also
$\Phi^{\circ r}$-periodic, the version of
Lemma~\ref{lem:lift} for $m$ coordinates gives an
irreducible $(F,\ldots,F)$-periodic curve $D$ with
dominant coordinate projections such that $C=\pi(D)$.
By the version of Proposition~\ref{cor:diagonal-param}
for $m$ coordinates and Theorem~\ref{thm:beaumont},
the curve $D$ is invariant under
$(F^{\circ n},\ldots,F^{\circ n})$. Hence
$$
 \Phi^{\circ rn}(C)
 =\pi\bigl((F^{\circ n},\ldots,F^{\circ n})(D)\bigr)
 =\pi(D)=C. \eqno{\Box} 
$$

\end{document}